\documentclass[12pt, twoside, leqno]{article}

\usepackage{amsmath,amsthm}
\usepackage{amssymb}

\usepackage{enumerate}

\newtheorem{thm}{Theorem}[section]
\newtheorem{cor}[thm]{Corollary}
\newtheorem{lem}[thm]{Lemma}
\newtheorem{alg}[thm]{Algorithm}
\newtheorem{prob}[thm]{Problem}

\theoremstyle{definition}
\newtheorem{defin}[thm]{Definition}
\newtheorem{rem}[thm]{Remark}

\newcommand{\Z}{\mathbb{Z}}
\newcommand{\N}{\mathbb{N}}

\newcommand{\PP}{\mathbb{P}}

\DeclareMathOperator{\sgn}{sgn}

\begin{document}


\baselineskip=17pt


\title{Computational aspects of rational residuosity}

\author{Markus Hittmeir\\
University of Salzburg\\
Mathematics Department}
\date{}

\maketitle


\renewcommand{\thefootnote}{}

\footnote{Address: Hellbrunnerstra{\ss}e 34, A-5020 Salzburg. E-Mail: markus.hittmeir@sbg.ac.at}



\footnote{The author is supported by the Austrian Science Fund (FWF): Project F5504-N26, which is a part of the Special Research Program "Quasi-Monte Carlo Methods: Theory and Applications".}

\renewcommand{\thefootnote}{\arabic{footnote}}
\setcounter{footnote}{0}


\vspace{-3cm}

\begin{abstract}
In this paper, we consider an extension of Jacobi's symbol, the so called rational $2^k$-th power residue symbol. In Section 3, we prove a novel generalization of Zolotarev's lemma. In Sections 4, 5 and 6, we show that several hard computational problems are polynomial-time reducible to computing these residue symbols, such as getting nontrivial information about factors of semiprime numbers. We also derive criteria concerning the Quadratic Residuosity Problem.
\end{abstract}

\section{Introduction}
The topic of this paper is the rational $2^k$-th power residue symbol.

\defin{Let $k\in\N_0$, $p\in\PP$ and $a\in\Z$ such that $p\nmid a$. We define the rational $2^k$-th power residue symbol as
\[
\left( \frac{a}{p} \right)_{2^k}:=\begin{cases} 1 &\mbox{if there is $x\in\Z$ such that } x^{2^k} \equiv a \mod p,\\
-1 & \mbox{otherwise.} \end{cases}
\]
If the symbol equals $1$, we say that $a$ is a \emph{$2^k$-th rational residue} modulo $p$. Let $n=p_1\cdots p_l$ be a product of not necessarily distinct primes such that $\gcd(n,a)=1$. We define
\[
\left( \frac{a}{n} \right)_{2^k}:=\left( \frac{a}{p_1} \right)_{2^k}\cdots \left( \frac{a}{p_l} \right)_{2^k}.
\]
We may also write this symbol as $(a|n)_{2^k}$.}\\

This symbol has been considered in several publications. Burde's and Scholz's reciprocity laws for the biquadratic case and their generalizations for the octic and for higher power cases are well-known. We refer the reader to \cite{L1}, to \cite[Theorem 3]{BCLS} and to \cite[Theorem 3.1]{BEK}. For a general overview on power residue symbols in number fields and other rational power residue symbols, we suggest the survey \cite{L2}.

The results mentioned above suppose that the prime modulus $p$ in the $2^k$-th residue symbol satisfies $p\equiv 1 \mod 2^k$. Usually, this assumption is part of the definition of the symbol itself. In Section 3, we demonstrate that there is particular interest in considering the symbol without such restriction. Theorem \ref{t3} extends the generalization of Zolotarev's lemma in \cite{BEKM} from the restricted notion to our general notion by dropping said assumption. Moreover, this extension is used in the proof of Theorem \ref{t5}, which is one of the main results in this paper and further generalizes Zolotarev's lemma for semiprime moduli. However, the remainder of the paper concerns the following computational problem.
{\prob{(Computation of the $2^k$-th residue symbol, $CRS(k)$)

Let $k\in\N$ be fixed. For any two coprime natural numbers $n,m$ satisfying $(m|p)_{2^{k-1}}=1$ for every prime factor $p$ of $n$, compute the value of the symbol 
\[
\left( \frac{m}{n} \right)_{2^k}.
\]
}\label{p1}}
Theorem \ref{t1} shows that this problem is easy to solve if the prime factorization of $n$ is known. But in general, this is not the case. Therefore, we are interested in methods to compute this symbol and solve the problem without knowing the factorization of $n$. We observe that there is a polynomial-time algorithm to solve $CRS(1)$, since the problem reduces to the computation of Jacobi's symbol, which can be done efficiently by using the quadratic reciprocity law. For $k>1$, the problem appears to be open and the currently known results and reciprocity laws do not apply. One of the reasons, the loss of generality due to the usual restriction to the case $p\equiv 1 \mod 2^k$, has already been mentioned. 

Since Zolotarev's lemma has been useful to prove quadratic reciprocity for Jacobi's symbol, Theorem \ref{t5} might also be of interest for attacking $CRS(k)$, $k>1$. However, this paper does not contain an efficient algorithm. On the contrary, our results indicate the computational hardness of $CRS$. We will prove polynomial-time reductions from three difficult problems to Problem \ref{p1}. In doing so, we also derive an interesting criterion for a case of the Quadratic Residuosity Problem in Corollary \ref{c3}.

To be concrete, Sections 4, 5 and 6 are respectively dedicated to the proofs of the following three results, where 'ERH' stands for 'Extended Riemann Hypothesis' and $\nu_2(n):=\max\{k\in\N_0:2^k\mid n\}$ denotes the $2$-adic valuation of $n\in\N$.
{\thm{Let $N\in\N$. There exists a probabilistic (and ERH-conditional deterministic) polynomial-time reduction from the problem of determining whether $N$ is the sum of two squares to $CRS(2)$.}}
{\thm{Let $N=pq$ be a semiprime number with unknown factors $p$ and $q$. There exists a probabilistic (and ERH-conditional deterministic) polynomial-time reduction from computing $v_p:=\nu_2(p-1)$ and $v_q:=\nu_2(q-1)$ to $CRS(k)$, where $k=1,2,...,\max\{v_p,v_q\}+1$.}}
{\thm{Let $N=pq$ be a semiprime number with unknown factors $p$ and $q$, where $\nu_2(p-1)\neq \nu_2(q-1)$. There exists a deterministic polynomial-time reduction from solving the Quadratic Residuosity Problem modulo $N$ to $CRS(\nu_2(N-1)+1)$.}}
\vspace{12pt}

It is an open question whether or not there is a polynomial-time reduction from the integer factorization problem to $CRS$.

\section{Basic properties}
The following theorem is a well-known generalization of Euler's Theorem about Legendre's symbol.

{\thm{Let $k\in\N$, $p\in\PP$ and $a\in\Z$ such that $p\nmid a$. Assume $(a|p)_{2^{k-1}}=1$. Then
\[
\left( \frac{a}{p} \right)_{2^k}\equiv a^{\frac{p-1}{\gcd(2^k,p-1)}}\mod p.
\]}\label{t1}}

\begin{proof}
Clearly, the statement is true for $p=2$. Let $p$ be an odd prime and $p-1=2^mu$ such that $u$ is odd. We consider two cases.

Case 1: $k\leq m$. Let $b$ be primitive modulo $p$ and $b^r\equiv a \mod p$ with $r\in\{1,...,p-1\}$. Then it is easy to see that $(a|p)_{2^k}=1$ if and only if $2^k\mid r$. Since $(a|p)_{2^{k-1}}=1$, we get $2^{k-1}\mid r$. We derive $(a|p)_{2^k}=(-1)^{r/2^{k-1}}$ and
\[
a^{\frac{p-1}{\gcd(2^k,p-1)}}\equiv a^{\frac{p-1}{2^k}}\equiv (b^r)^{\frac{p-1}{2^k}}\equiv (b^{\frac{p-1}{2}})^{\frac{r}{2^{k-1}}}\equiv (-1)^{\frac{r}{2^{k-1}}} \mod p.
\]
Case 2: $k>m$. From $(a|p)_{2^{k-1}}=1$ it easily follows that $(a|p)_{2^m}=1$. In Case $1$ we showed that this implies 
\[
a^{\frac{p-1}{\gcd(2^k,p-1)}}\equiv a^{\frac{p-1}{2^m}}\equiv 1\mod p.
\] 
Therefore, we have to show that $(a|p)_{2^k}=1$. We know that there is $x\in\Z$ such that $x^{2^m}\equiv a \mod p$. Define $t:=(2^{k-m})^{-1}\pmod u$ and $y:=x^t$. Since $u$ divides $t2^{k-m}-1$, $p-1$ divides 
\[
2^m(t2^{k-m}-1)=t2^k-2^m
\]
and we derive $x^{t2^k-2^m}\equiv 1 \mod p$, hence
\[
y^{2^k}\equiv x^{2^m}\equiv a \mod p,
\]
which we wanted to show.
\end{proof}

{\lem{Let $n\in\N$ and $a,b\in\Z$ such that $\gcd(n,ab)=1$. Then the following holds:
\begin{enumerate}
\item{If $a\equiv b\mod n$, then $(a|n)_{2^k}=(b|n)_{2^k}$ for every $k\in\N_0$.}
\item{Let $k\in\N$ and $(a|p)_{2^{k-1}}=(b|p)_{2^{k-1}}=1$ for every prime factor $p$ of $n$. Then
\[
\left( \frac{ab}{n} \right)_{2^k}=\bigg( \frac{a}{n} \bigg)_{2^k}\cdot\left( \frac{b}{n} \right)_{2^k}.
\]}
\end{enumerate}}}

\begin{proof}
1: By definition, $(a|p)_{2^k}=(b|p)_{2^k}$ holds for every prime factor $p$ of $n$. This implies the claim.

2: We prove this equality by induction. For $k=1$, the claim follows from the multiplicativity of Jacobi's symbol. Assume that the statement holds for $k-1$, then $(ab|p)_{2^{k-1}}=(a|p)_{2^{k-1}}\cdot(b|p)_{2^{k-1}}=1$ for every prime factor $p$ of $n$. We apply Theorem \ref{t1} and deduce
\[
\left( \frac{ab}{p} \right)_{2^k}\equiv (ab)^{\frac{p-1}{\gcd(2^k,p-1)}}\equiv \left( \frac{a}{p} \right)_{2^k}\cdot\left( \frac{b}{p} \right)_{2^k}\mod p
\]
for every prime factor $p$ of $n$, which implies the statement.
\end{proof}

\rem{For $k\in\N_0$ and $n\in\N$, the second property in the preceding lemma yields that the $2^k$-th rational residues modulo $n$ build a subgroup of $\Z_n^*$. We will denote this subgroup by $\Z_{n,2^k}^*.$} 

{\lem{Let $p\in\PP$ and $a\in\Z$ such that $p\nmid a$. Let $m:=\nu_2(p-1)$ be the $2$-adic valuation of $p-1$, then
\[
\left( \frac{a}{p} \right)_{2^k}=\left( \frac{a}{p} \right)_{2^m}
\]
for every $k\geq m$.\label{l1}}}

\begin{proof}
We will show that $\Z_{p,2^k}^{*}=\Z_{p,2^m}^{*}$ holds for every $k\geq m$. The stated equality is an immediate consequence of this fact. Let $k\geq m$ be arbitrary. If $a\in\Z_{p,2^k}^{*}$, then $(a|p)_{2^k}=1$. By definition, this obviously implies $(a|p)_{2^m}=1$. We get $a\in\Z_{p,2^m}^{*}$ and therefore $\Z_{p,2^k}^{*}\subseteq \Z_{p,2^m}^{*}$. 

We will prove $\Z_{p,2^m}^{*}\subseteq \Z_{p,2^k}^{*}$ by induction. For $k=m$, there is nothing to show. Now let $a\in\Z_{p,2^m}^{*}$ and assume that the statement holds for the exponent $k-1$. We derive $a\in\Z_{p,2^{k-1}}^*$, hence $(a|p)_{2^{k-1}}=1$. Theorem \ref{t1} yields
\[
\left( \frac{a}{p} \right)_{2^k}\equiv a^{\frac{p-1}{\gcd(2^k,p-1)}}\equiv a^{\frac{p-1}{2^m}}\mod p.
\]
Since $(a|p)_{2^m}=1$ implies $(a|p)_{2^{m-1}}=1$, Theorem \ref{t1} also yields 
\[
a^{\frac{p-1}{2^m}}\equiv \left( \frac{a}{p} \right)_{2^m}=1 \mod p.
\]
We deduce $(a|p)_{2^k}=1$ and $a\in\Z_{p,2^k}^{*}$, which we wanted to show.
\end{proof}

{\cor{Let $p\in\PP$ and $m:=\nu_2(p-1)$. Then $\Z_{p,2^k}^{*}=\Z_{p,2^m}^{*}$ for every $k\geq m$.}\label{c1}}\\

To prove the results concerning the complexity of computing the $2^k$-th residue symbol in the Sections 4,5 and 6, we will use the following lemma on a regular basis.

{\lem{Let $k\in\N$, $p\in\PP$ and $a\in\Z$ such that $p\nmid a$. Then
\[
\left( \frac{a^{2^{k-1}}}{p} \right)_{2^k}=\begin{cases}\left( \frac{a}{p} \right)_2 &\mbox{if } k\leq \nu_2(p-1), \\
1 & \mbox{otherwise.}\end{cases}
\]}\label{l3}}

\begin{proof}
We observe that $\left( \frac{a^{2^{k-1}}}{p} \right)_{2^{k-1}}=1$. Therefore, Theorem \ref{t1} implies
\[
\left( \frac{a^{2^{k-1}}}{p} \right)_{2^k}\equiv a^{\frac{2^{k-1}(p-1)}{\gcd(2^k,p-1)}}\equiv \begin{cases}\left( \frac{a}{p} \right)_2 \mod p &\mbox{if } k\leq \nu_2(p-1), \\
1 \mod p & \mbox{otherwise.}\end{cases}
\]
Since the symbols only take on the values $1$ and $-1$, this already proves the claim.
\end{proof}

\section{Generalizations of Zolotarev's lemma}

In the following, we will describe a way to express the $2^k$-th residue symbol via the signum of a certain permutation on $\Z_{p}$. The statement for $k=1$ is equivalent to Zolotarev's lemma, which has been used by Y. I. Zolotarev in his proof of quadratic reciprocity in 1872. The following theorem has already been proven in \cite{BEKM} for the restricted notion of the $2^k$-th residue symbol. We now show that we may drop the assumption $p\equiv 1 \mod 2^k$ and that, hence, the claim also holds for the general notion.

{\thm{Let $k\in\N$, $p\in\PP$ and $a\in\Z$ such that $p\nmid a$. Assume $(a|p)_{2^{k-1}}=1$. Then
\[
\left( \frac{a}{p} \right)_{2^k}=\sgn\left(\phi_a|_{\Z_{p,2^{k-1}}^{*}}\right),
\]
where $\phi_a$ is the permutation on $\Z_{p}$ given by 
\[
x \pmod {p}\mapsto ax \pmod {p}.
\]}\label{t3}}

\begin{proof}
It is easy to see that the statement is true for $p=2$. Let $p$ be an odd prime. The case $p\equiv 1 \mod 2^k$ has been proven in \cite[Theorem 1]{BEKM}. We assume $p\not\equiv 1 \mod 2^k$. Let $m:=\nu_2(p-1)$ be the $2$-adic valuation of $p-1$ and note that $k>m$. From the assumption and Lemma \ref{l1} it follows that
\[
\left( \frac{a}{p} \right)_{2^k}=\left( \frac{a}{p} \right)_{2^{k-1}}=1.
\]
Therefore, there exists $z\in\Z_{p}^*$ such that $z^{2^k}\equiv a \mod p$. Obviously, the element $b:=z^{2^{k-1}}$ is in $\Z_{p,2^{k-1}}^*$ and $b^2\equiv a \mod p$ holds. Now $\sigma_a:=\phi_a|_{\Z_{p,2^{k-1}}^{*}}$ and $\sigma_b:=\phi_b|_{\Z_{p,2^{k-1}}^{*}}$ are both permutations on $\Z_{p,2^{k-1}}^*$ and we have
\[
\sigma_a(x)=ax\pmod {p}=b^2x \pmod {p}=b(bx) \pmod {p}
\]
for every $x\in\Z_{p,2^{k-1}}^{*}$, hence $\sigma_a=\sigma_b\circ\sigma_b$. By using the multiplicativity of the signum of compositions of permutations, we conclude that
\[
\sgn \sigma_a=\sgn (\sigma_b\circ \sigma_b)=(\sgn \sigma_b)^2=1,
\]
and the statement follows.
\end{proof}

\rem{Let $k\in\N_0$ and $n\in\N$. We will continue with proving a further generalization of Zolotarev's lemma for semiprime moduli. In addition to the already defined subgroups $\Z_{n,2^k}^*$ of $\Z_n^*$, we will also consider
\[
\Z_{n,2^k}:=\{r\in\Z_n\mid \exists x \in\Z :x^{2^k}\equiv r\mod n\}.
\]
}

We will make use of an elementary fact about permutations on cartesian products. The reader may also find Lemma \ref{l5} in \cite{DS}.

{\lem{
Let $X_1,...,X_l$ be finite sets and $\sigma_i\in Sym(X_i)$, $i=1,...,l$. For $X:=X_1\times X_2\times...\times X_l$, define the permutation $\sigma\in Sym(X)$ by applying the $\sigma_i$ componentwise to the $l$-tupels in X. Then 
\[
\sgn(\sigma)=\prod_{i=1}^l \sgn(\sigma_i)^{\gamma_i},
\]
 where $\gamma_i:=|X|/|X_i|$, $i=1,...,l$.
}\label{l5}}

\begin{proof}
Let $\hat{\sigma_i}$ be obtained by applying $\sigma_i$ to the $i$-th coordinate and the identity to the other coordinates. Then $\hat{\sigma_i}$ consists of $\gamma_i$ copies of $\sigma_i$, which implies $\sgn(\hat{\sigma_i})=\sgn(\sigma_i)^{\gamma_i}$. As a consequence of $\sigma=\hat{\sigma_1}\circ \hat{\sigma_2}\circ...\circ \hat{\sigma_l}$, the statement follows.
\end{proof}

For $n,b\in\N$, let $\nu_b(n):=\max\{k\in\N_0:b^k\mid n\}$ be the $b$-adic valuation of $n$. We will need the following result.

{\lem{Let $b\in\N$ and $N=pq$, where $p,q\in\PP$. W.l.o.g., we suppose $\nu_b(p-1)\leq \nu_b(q-1)$. Then:
\[
\nu_b(p-1)\begin{cases} =\nu_b(N-1) &\mbox{if } \nu_b(p-1)<\nu_b(q-1), \\
\leq \nu_b(N-1) & \mbox{otherwise.}\end{cases} 
\]
}\label{l2}}

\begin{proof}
Let $\nu_b(p-1)<\nu_b(q-1)$ and $k\in\N_0$ such that $k\leq \nu_b(p-1)$. We clearly have $b^k\mid p-1$ and $b^k\mid q-1$. Since $k< \nu_b(q-1)$, we also derive $(q-1)/b^k\equiv 0 \mod b$ and
\[
\frac{p-1}{b^k}\equiv \frac{p-1}{b^k} +\frac{q-1}{b^k} \mod b.
\]
Furthermore, $\nu_b(q-1)$ must be larger than $0$, which implies $q\equiv 1 \mod b$. We obtain
\[
\frac{p-1}{b^k} +\frac{q-1}{b^k} \equiv q \frac{p-1}{b^k} + \frac{q-1}{b^k}= \frac{N-1}{b^k} \mod b.
\]
We conclude that $(N-1)/b^k \equiv (p-1)/b^k \mod b$ for every $k\leq \nu_b(p-1)$ and, therefore, $\nu_b(p-1)=\nu_b(N-1)$.

Now let $\nu_b(p-1)=\nu_b(q-1)=:k$. We have
\[
N=pq\equiv 1\cdot 1 =1 \mod b^k,
\]
hence $b^k \mid N-1$, which implies that $k\leq \nu_b(N-1)$.
\end{proof}

\rem{If we consider $b=2$, then it is easy to prove that in the case $\nu_2(p-1)=\nu_2(q-1)=:k$, we have $N\equiv 1 \mod 2^{k+1}$ and therefore $\nu_2(p-1)<\nu_2(N-1)$.\label{r1}}

{\thm{
Let $k\in\N$ and $N=pq$, where $p,q$ are odd primes and $p\neq q$. Let $m\in\N$ be coprime to $N$ such that $(m|p)_{2^{k-1}}=(m|q)_{2^{k-1}}=1$. Then
\[
\left( \frac{m}{N} \right)_{2^k}=\begin{cases} \sgn\left(\phi_m|_{\Z_{N,2^{k-1}}}\right) &\mbox{if } N\equiv 1 \mod 2^k, \\
\sgn\left(\phi_m|_{\Z_{N,2^{k-1}}^*}\right) & \mbox{if } N\not\equiv 1 \mod 2^k,\end{cases}
\]
where $\phi_m$ is the permutation on $\Z_{N}$ given by 
\[
x \pmod {N}\mapsto mx \pmod {N}.
\]
}\label{t5}}

\begin{proof}
W.l.o.g., we may assume that $\nu_2(p-1)\leq \nu_2(q-1)$.

The case $N\not\equiv 1 \mod 2^k$: 
Define $\sigma_1^*:=\phi_m|_{\Z_{p,2^{k-1}}^*}$ and $\sigma_2^*:=\phi_m|_{\Z_{q,2^{k-1}}^*}$. Let $\sigma^*$ be the permutation on $X:=\Z_{p,2^{k-1}}^*\times \Z_{q,2^{k-1}}^*$ defined by applying $\sigma_1^*$ and $\sigma_2^*$ componentwise to pairs in $X$. If $\psi$ is the bijection from $\Z_{N,2^{k-1}}^*$ to $X$ given by the Chinese Remainder Theorem, then one easily observes that
\[
\phi_m|_{\Z_{N,2^{k-1}}^*}=\psi^{-1}\circ\sigma^*\circ\psi,
\]
which implies $\sgn\left(\phi_m|_{\Z_{N,2^{k-1}}^*}\right)=\sgn(\sigma^*)$. 

Note that we have $k>\nu_2(N-1)$. From Lemma \ref{l2} it easily follows that one of the following two cases hold:
\begin{enumerate}
\item{$p\not\equiv 1\mod 2^k \wedge q\equiv 1 \mod 2^k$: We are going to apply Lemma \ref{l5}. Setting $\gamma_1^*:=|\Z_{q,2^{k-1}}|$ and $\gamma_2^*:=|\Z_{p,2^{k-1}}|$, we deduce
\[
\sgn(\sigma^*)=\sgn(\sigma_1^*)^{\gamma_1^*}\cdot \sgn(\sigma_2^*)^{\gamma_2^*}.
\]
Since $\gamma_1^*=(q-1)/2^{k-1}\equiv 0 \mod 2$ and $\gamma_2^*\equiv 1 \mod 2$, we get $\sgn(\sigma^*)=\sgn(\sigma_2^*)$. By using Theorem \ref{t1} and Theorem \ref{t3},we conclude that 
\[
\left( \frac{m}{N} \right)_{2^k}=\left( \frac{m}{p} \right)_{2^k}\cdot \left( \frac{m}{q} \right)_{2^k}=\left( \frac{m}{q} \right)_{2^k}=\sgn(\sigma_2^*).
\]}
\item{$p\not\equiv 1 \mod 2^k \wedge q\not\equiv 1\mod 2^k$: Again, we apply Lemma \ref{l5}. We derive
\begin{align*}
\sgn(\sigma^*)=\sgn(\sigma_1^*)^{\gamma_1^*}\cdot \sgn(\sigma_2^*)^{\gamma_2^*}&=\sgn(\sigma_1^*)\cdot\sgn(\sigma_2^*)\\
&=\left( \frac{m}{p} \right)_{2^k}\cdot \left( \frac{m}{q} \right)_{2^k}=\left( \frac{m}{N} \right)_{2^k},
\end{align*}
where we have used $\gamma_1^*\equiv\gamma_2^*\equiv 1 \mod 2$ and Theorem \ref{t3}.
}
\end{enumerate}

The case $N\equiv 1 \mod 2^k$:
Analogous to the previous case, we define $\sigma_1:=\phi_m|_{\Z_{p,2^{k-1}}}$, $\sigma_2:=\phi_m|_{\Z_{q,2^{k-1}}}$ and $\sigma$ as the permutation on the set $X:=\Z_{p,2^{k-1}}\times \Z_{q,2^{k-1}}$ which is given by applying $\sigma_1$ and $\sigma_2$ componentwise to pairs in $X$. Again, it is easy to see that the Chinese Remainder Theorem yields a bijection from $\Z_{N,2^{k-1}}$ to $X$ and that we may derive $\sgn\left(\phi_m|_{\Z_{N,2^{k-1}}}\right)=\sgn(\sigma)$. 

Now we have $k\leq \nu_2(N-1)$. Assuming $2^k\nmid p-1$ and $2^k\mid q-1$, it is easy to deduce a contradiction by using Lemma \ref{l2}. Hence, we are left with the following two cases:
\begin{enumerate}
\item{$p\equiv q\equiv 1 \mod 2^k$:
Applying Lemma \ref{l5} with $\gamma_1:=|\Z_{q,2^{k-1}}|$ and with $\gamma_2:=|\Z_{p,2^{k-1}}|$, we get
\[
\sgn(\sigma)=\sgn(\sigma_1)^{\gamma_1}\cdot\sgn(\sigma_2)^{\gamma_2}=\sgn(\sigma_1)\cdot\sgn(\sigma_2),
\]
where we have used $\gamma_1=(q-1)/2^{k-1}+1$, $\gamma_2=(p-1)/2^{k-1}+1$ and, hence, $\gamma_1\equiv \gamma_2\equiv 1 \mod 2$. Theorem \ref{t3} yields
\[
\left( \frac{m}{N} \right)_{2^k}=\left( \frac{m}{p} \right)_{2^k}\cdot \left( \frac{m}{q} \right)_{2^k}=\sgn(\sigma_1^*)\cdot\sgn(\sigma_2^*).
\]
Since one easily observes that $\sgn(\sigma_i^*)=\sgn(\sigma_i)$ for $i=1,2$, this proves the statement.
}
\item{$p\not\equiv 1 \mod 2^k \wedge q\not\equiv 1\mod 2^k$:
Lemma \ref{l5} now yields 
\[
\sgn(\sigma)=\sgn(\sigma_1)^{\gamma_1}\cdot\sgn(\sigma_2)^{\gamma_2}=1,
\]
since $\gamma_1\equiv \gamma_2\equiv 0 \mod 2$. By Theorem \ref{t1}, we also deduce
\[
\left( \frac{m}{N} \right)_{2^k}=1.
\]
}
\end{enumerate}
This concludes the proof.
\end{proof}

\rem{There is a version of Zolotarev's lemma for Jacobi's symbol, stating that $(m|n)_2=\sgn(\phi_m)$ for coprime $m,n\in\N$, where $n>1$ is odd and $\phi_m$ is considered as permutation on $\Z_n$. For semiprime moduli $N$, Theorem \ref{t5} extends this result. However, it does not hold for any composite moduli. Consider, for example, a product of three primes $n=pqr$ with $p\equiv 3\mod 4$ and $r,q\equiv 1 \mod 4$. Using Lemma \ref{l5}, it is an easy exercise to show that $\sgn\left(\phi_m|_{\Z_{n,2}}\right)=\sgn\left(\phi_m|_{\Z_{n,2}^*}\right)=1$ for any $m$ satisfying our assumptions; yielding a counterexample for every $m$ with $(m|n)_4=-1$.

\section{A solvability criterion for $N=X^2+Y^2$}
Let $N$ be any natural number. In this section, we discuss how an efficient algorithm for computing the biquadratic rational power residue symbol could be used to efficiently determine whether the diophantine equation $N=X^2+Y^2$ has a solution or not. The following theorem is a well-known statement.

{\thm{(P. Fermat)

Let $N\in\N$. The diophantine equation $N=X^2+Y^2$ has a solution if and only if every prime factor $p$ of $N$ with $p\equiv 3 \mod 4$ occurs to an even power in the prime factorization of $N$.}}\\

We will use this statement to prove a solvability criterion for this equation via $4$-th residue symbols. We first consider the following lemma, which yields a necessary condition for solvability.

{\lem{Let $N\in\N$ and $a\in\Z$ with $\gcd(N,a)=1$. If $N=X^2+Y^2$ is solvable over the integers, then 
\[
\left( \frac{a^2}{N} \right)_{4}=\left( \frac{a}{N} \right)_{2}.
\]}\label{l4}}
\begin{proof}
If $N=X^2+Y^2$ is solvable over the integers, Fermat's result allows us to denote 
\[
N=2^{e_0}p_1^{2e_1}\cdots p_{v-1}^{2e_{v-1}}p_v^{e_v}\cdots p_r^{e_r}, 
\]
where $p_1,...,p_{v-1}$ are precisely those prime factors with $p_i\equiv 3 \mod 4$. We derive
\begin{align*}
\left( \frac{a}{N} \right)_{2}=&\left( \frac{a}{2} \right)_{2}^{e_0}\cdot \left( \frac{a}{p_1^{e_1}\cdots p_{v-1}^{e_{v-1}}} \right)_{2}^2 \cdot\left( \frac{a}{p_v^{e_v}\cdots p_r^{e_r}} \right)_{2}\\
=&\left( \frac{a}{p_v^{e_v}\cdots p_r^{e_r}} \right)_{2}=\left( \frac{a}{p_v} \right)_{2}^{e_v}\cdots \left( \frac{a}{p_r} \right)_{2}^{e_r}\\
=&\left( \frac{a^2}{p_v} \right)_{4}^{e_v}\cdots \left( \frac{a^2}{p_r} \right)_{4}^{e_r}=\left( \frac{a^2}{p_v^{e_v}\cdots p_r^{e_r}} \right)_{4}.
\end{align*}
where we have used Lemma \ref{l3} and $p_i\equiv 1 \mod 4$ for $i=v,...,r$. Note that one easily observes
\[
\left( \frac{a^2}{2^{e_0}p_1^{2e_1}\cdots p_{v-1}^{2e_{v-1}}} \right)_{4}=1,
\]
which implies the claim.
\end{proof}

Our criterion relies on the Extended Riemann Hypothesis. In the proof, we will need the following result about least quadratic nonresidues, which has been proven in \cite[Theorem 6.35]{W}.

{\thm{(S. Wedeniwski, ERH)

Assume that the Extended Riemann Hypothesis is correct. Let $m$ be an odd positive integer greater than $1$ such that $m$ is not a perfect square, and 
$
x:=\min\{k\in\N : \left( \frac{k}{m} \right)_{2}\neq 1\}.
$
Then 
\[
x<\frac{3}{2}(\log m)^2-\frac{44}{5}\log m+13.
\]
}}

\rem{Currently, we do not know an unconditional logarithmic bound for the least quadratic nonresidue.}

{\thm{(ERH) Let $N\in\N$, set
\[
M:=\Big{\{}p\in\PP:p<\frac{3}{2}(\log N)^2-\frac{44}{5}\log N+13\Big{\}}
\]
and assume that $p\nmid N$ for every $p\in M$. If the Extended Riemann Hypothesis is correct, then
\[
N=X^2+Y^2\textit{ is solvable in $\Z$ } \Longleftrightarrow \forall a \in M: \left( \frac{a^2}{N} \right)_{4}=\left( \frac{a}{N} \right)_{2}.
\]}\label{t2}}

\begin{proof}
Lemma \ref{l4} shows that the right statement follows from the left one. Now we assume to the contrary that $N=X^2+Y^2$ has no integer solutions. Since $2\in M$, the assumptions imply that $N$ is odd. According to Fermat's result, we denote
\[
N=p_1^{e_1}\cdots p_{v-1}^{e_{v-1}}p_v^{2e_v}\cdots p_{\mu-1}^{2e_{\mu-1}}p_{\mu}^{e_{\mu}}\cdots p_r^{e_r},
\]
where $p_i\equiv 3 \mod 4$ for $i=1,...,\mu-1$ and $p_i\equiv 1 \mod 4$ for $i=\mu,...,r$. Since the diophantine equation is not solvable, there is at least one prime factor $p$ of $N$ with $p\equiv 3 \mod 4$ that occurs to an odd power in the prime factorization. Hence, we assume $e_i\equiv 1\mod 2$ for $i=1,...,v-1$, $v\geq 2$. Let $a\in\Z$ with $\gcd(N,a)=1$ be arbitrary. We consider
\begin{align*}
\left( \frac{a}{N} \right)_{2}=&\left( \frac{a}{p_1} \right)_{2}^{e_1}\cdots \left( \frac{a}{p_{v-1}}\right)_{2}^{e_{v-1}}  \cdot\left( \frac{a}{p_v^{2e_v}\cdots p_{\mu-1}^{2e_{\mu-1}}p_\mu^{e_\mu}\cdots p_{r}^{e_{r}}} \right)_{2}\\
=&\left( \frac{a}{p_1\cdots p_{v-1}} \right)_{2}\cdot \left( \frac{a}{p_v^{2e_v}\cdots p_{\mu-1}^{2e_{\mu-1}} p_\mu^{e_\mu}\cdots p_{r}^{e_{r}}} \right)_{2}.
\end{align*}
Here, we have used that $e_i\equiv 1 \mod 2$ for $i=1,...,v-1$. According to Fermat's result and Lemma \ref{l4}, we derive
\[
\left( \frac{a}{N} \right)_{2}=\left( \frac{a}{p_1\cdots p_{v-1}} \right)_{2}\cdot \left( \frac{a^2}{p_v^{2e_v}\cdots p_{\mu-1}^{2e_{\mu-1}}p_\mu^{e_\mu}\cdots p_{r}^{e_{r}}} \right)_{4}.
\]
It is easy to show that, modulo any number $m$, the smallest element $x$ satisfying $(x|m)_2=-1$ must be prime. Hence, from Wedeniwski's result and our assumption that $\gcd(p,N)=1$ for $p\in M$ it follows that there is some $a_0\in M$ such that 
\[
(a_0|p_1\cdots p_{v-1})_2=-1.
\]
We derive
\[
\left( \frac{a_0}{N} \right)_{2}=- \left( \frac{a_0^2}{p_v^{2e_v}\cdots p_{\mu-1}^{2e_{\mu-1}}p_\mu^{e_\mu}\cdots p_{r}^{e_{r}}} \right)_{4}\neq \left( \frac{a_0^2}{p_v^{2e_v}\cdots p_{\mu-1}^{2e_{\mu-1}}p_\mu^{e_\mu}\cdots p_{r}^{e_{r}}} \right)_{4}.
\]
But for $i=1,...,v-1$, Lemma \ref{l3} yields $(a_0^2|p_i)_4=1$, which implies 
\[
\left( \frac{a_0^2}{p_1^{e_1}\cdots p_{v-1}^{e_{v-1}}} \right)_{4}=1.
\]
We conclude that $(a_0|N)_2\neq (a_0^2|N)_4$, which we wanted to show.
\end{proof}

\rem{\begin{enumerate}
\item{One observes that we can easily deal with the assumption of Theorem \ref{t2}. We simply perform trial division up to Wedeniwski's bound. If any prime factor of $N$ is found, we test if it satisfies the conditions of Fermat's result. If this is not the case, the diophantine equation is not solvable. If the conditions are satisfied, we remove this prime factor from $N$ and the question about the solvability of the equation corresponding to the original number $N$ reduces to the question about the solvability of the equation corresponding to the resulting number. Obviously, this procedure can be performed in polynomial time.}
\item{Theorem \ref{t2} suggests a deterministic algorithm, but the proof of its correctness relies on the ERH. We want to point out that we could also use a probabilistic argument based on random choices of $a$. Note that half of all the numbers $a\in\Z_N^*$ satisfy $(a|p_1\cdots p_{v-1})_2=-1$. Assume that the equation $N=X^2+Y^2$ is not solvable in $\Z$. Then after $k$ choices for $a$, the probability of not having found any element which does not satisfy our criterion is $2^{-k}$.}
\end{enumerate}}

\section{The first bits of the factors of a semiprime}
Let $N=pq$ be a natural number, where $p,q$ are large, unknown and distinct primes. The security of the RSA-cryptosystem relies on the difficulty of the problem to compute $p$ and $q$ if only $N$ is known. In this section, we show how an efficient algorithm for the computation of the rational $2^k$-th power residue symbol could be used to efficiently compute the $2$-adic valuations $\nu_2(p-1)$ and $\nu_2(q-1)$ and, therefore, the first few bits of the factors of $N$. 

{\alg{
Let $N=pq$ be odd, $p,q\in\PP$. W.l.o.g., we suppose that $\nu_2(p-1)\leq \nu_2(q-1)$. Take the following steps to compute $\nu_2(p-1)$ and $\nu_2(q-1)$.
\begin{enumerate}
\item{Compute $v:=\nu_2(N-1)$. Go to Step 2.}
\item{Find $a\in\Z$ such that $(a|N)_2=-1$ via random search. Go to Step 3.}
\item{For $i=1,...,v$, compute $s_i:=(a^{2^{i-1}}|N)_{2^i}$. If there exists a minimal $j\in\{1,...,v\}$ with $s_j=1$, set 
\[
\nu_2(p-1)=\nu_2(q-1)=j-1
\]
and stop. If $s_i=-1$ for all $i\in\{1,...,v\}$, set $\nu_2(p-1)=v$ and go to Step 4.}
\item{Find $b\in\Z$ such that $(b^{2^v}|N)_{2^{v+1}}=-1$ via random search. Go to Step 5.}
\item{For $i\geq v+1$, compute $s_i:=(b^{2^{i-1}}|N)_{2^{i}}$. Do until a minimal $j\geq v+1$ is found with $s_j=1$. Set $\nu_2(q-1)=j-1$ and stop.}
\end{enumerate}\label{a1}}

\begin{proof}[Proof of Correctness]
The case $k:=\nu_2(p-1)=\nu_2(q-1)$:
According to Lemma \ref{l2} and Remark \ref{r1}, $k<v$ holds. We consider $s_i$ in Step 3. Lemma \ref{l3} yields that
\[
s_i=\left( \frac{a^{2^{i-1}}}{p} \right)_{2^i}\cdot \left( \frac{a^{2^{i-1}}}{q} \right)_{2^i}=\begin{cases}\left( \frac{a}{N} \right)_2 &\mbox{if } i\leq k, \\
1 & \mbox{otherwise.}\end{cases}
\]
Therefore, there exists a minimal $j\in\{1,...,v\}$ such that $s_j=(a|N)_2=1$, namely $j=k+1$. We conclude that, in this case, the algorithm correctly computes $k$ in Step 3.

The case $\nu_2(p-1)<\nu_2(q-1)$: 
It follows from Lemma \ref{l2} that we have $v=\nu_2(p-1)<\nu_2(q-1)$. We apply Lemma \ref{l3}, and for all $i\in\{1,...,v\}$, we derive
\[
s_i=\left( \frac{a}{N} \right)_2=-1.
\]
Hence, the algorithm correctly sets $\nu_2(p-1)=v$ in Step 3. Note that since $v=\nu_2(p-1)$ and $v+1\leq\nu_2(q-1)$, we obtain
\[
\left(\frac{b^{2^v}}{N} \right)_{2^{v+1}}=\left(\frac{b^{2^v}}{p} \right)_{2^{v+1}}\cdot \left(\frac{b^{2^v}}{q} \right)_{2^{v+1}}=\left(\frac{b}{q} \right)_2,
\]
where we have used Lemma \ref{l3} again. As a consequence, the algorithm actually finds $b\in\Z$ which is a quadratic nonresidue modulo $q$ in Step 4. Considering the values of $s_i$ in Step 5 in a similar manner, we deduce
\[
s_i=\begin{cases}\left( \frac{b}{q} \right)_2 &\mbox{if } v+1\leq i\leq \nu_2(q-1), \\
1 & \mbox{otherwise.}\end{cases}
\]
Now the minimal $j\geq v+1$ with $s_j=1$ obviously is $j=\nu_2(q-1)+1$.
\end{proof}

\rem{(Running Time of Algorithm \ref{a1})
\begin{enumerate}
\item{Algorithm \ref{a1} uses probabilistic methods in Step 2 and Step 4 and, therefore, is not deterministic. However, we want to point out that the probabibility of not finding a suitable element in these steps equals $2^{-k}$ after $k$ trials. This is due to our observation about the symbol in Step 4 in the correctness proof and due to elementary results on Legendre's and Jacobi's symbol.}
\item{It is easy to see that the running time of the algorithm mainly depends on the complexity of computing the rational power residue symbol.The other computations can be done in polynomial time and all the loops finish after at most $O(\log N)$ runs; either provably or at least with high probability.} 
\item{Relying on the Extended Riemann Hypothesis, we may formulate a deterministic version of the algorithm by using Wedeniwski's bound discussed in Section 4.}
\end{enumerate}
}

\rem{
It is easy to see that, if we know the valuations $v_p:=\nu_2(p-1)$ and $v_q:=\nu_2(q-1)$, we also have found a fast way to compute $p \pmod {2^m}$ and $q \pmod {2^m}$ for $m:=\max \{v_p,v_q\}+1$. 
}

\section{The Quadratic Residuosity Problem}
Several cryptographic methods rely on the computational hardness of the following problem.

{\prob{(Quadratic Residuosity Problem, QRP)

Let $N=pq$ with unknown, different primes $p$ and $q$ and let $a\in\Z$ such that $(a|N)_2=1$. Determine whether $a\in\Z_{N,2}^*$ or not.}}\\

QRP is about deciding whether $a$ is a quadratic residue modulo $N$ or not. The following theorem is an immediate consequence of our previously shown statements. It demonstrates how a certain case of this problem could be solved if we were able to efficiently compute rational $2^k$-th power residue symbols.

{\thm{Let $N=pq$ with $p,q\in\PP$ such that $\nu_2(p-1)\neq\nu_2(q-1)$. Let $a\in\Z$ with $(a|N)_2=1$. For $v:=\nu_2(N-1)$, it holds that
\[
a \in\Z_{N,2}^* \Longleftrightarrow \left( \frac{a^{2^{v}}}{N} \right)_{2^{v+1}}=1.
\]}\label{t4}}

\begin{proof} We assume $\nu_2(p-1)<\nu_2(q-1)$ w.l.o.g. Lemma \ref{l2} yields that $v=\nu_2(p-1)<\nu_2(q-1)$. Using Lemma \ref{l3}, we derive
\[
\left( \frac{a^{2^{v}}}{N} \right)_{2^{v+1}}=\left( \frac{a}{q} \right)_2.
\]
Since it is easy to show that $a$ is a quadratic residue modulo $N$ if and only if $(a|q)_2=1$, the claim is proven.
\end{proof}

{\cor{Let $N=pq$ be semiprime with distinct factors $p$ and $q$ such that $N\equiv 3 \mod 4$. Let $a\in\Z$ with $(a|N)_2=1$. Then it holds that
\[
a \in\Z_{N,2}^* \Longleftrightarrow  \left( \frac{a^{2}}{N} \right)_{4}=1.
\]}\label{c2}}

\begin{proof}
One easily observes that the pairs $(3,1)$ and $(1,3)$ are the only solutions to the congruence $N\equiv XY\mod 4$. Therefore, one of the prime factors of $N$ is congruent to $3$ modulo $4$, whereas the other one is congruent to $1$ modulo $4$. This implies $\nu_2(p-1)\neq \nu_2(q-1)$. Since $\nu_2(N-1)=1$, we may apply Theorem \ref{t4} to conclude the proof.
\end{proof}

\rem{Corollary \ref{c2} yields that for semiprimes $N$ with $N\equiv 3 \mod 4$, there is a deterministic polynomial-time reduction from solving the $QRP$ modulo $N$ to $CRS(2)$.}\\

We are now able to deduce the following criterion for quadratic residuosity modulo semiprimes $N$ satisfying $N\equiv 3\mod 4$.

{\cor{Let $N=pq$ be semiprime with distinct factors $p$ and $q$ such that $N\equiv 3 \mod 4$. Let $a\in\Z$ with $(a|N)_2=1$. Then it holds that
\[
a \in\Z_{N,2}^* \Longleftrightarrow  \sgn\left(\phi_{a^2}|_{Z_{N,2}^*}\right)=1,
\]
where $\phi_{a^2}$ is the permutation on $\Z_{N}^*$ given by 
\[
x \pmod {N}\mapsto a^2x \pmod {N}.
\]
}\label{c3}}

\begin{proof}
This statement is a consequence of Corollary \ref{c2} applied together with Theorem \ref{t5}.
\end{proof}

\end{document}